\documentstyle[12pt,leqno]{article}

\begin{document}
\title{Extension Dimension and Refinable Maps}
\author{Alex Chigogidze and Vesko Valov\\[3pt] }

\date{}

\maketitle

{\bf Abstract.}
Extension dimension is characterized in terms of $\omega$-maps.  
We apply
this result to prove that extension dimension is preserved by refinable
maps between metrizable spaces. It is also shown that refinable maps
preserve some infinite-dimensional properties.

\bigskip\noindent
1991 MS Classification: 55M10; 54C10.infinite-dimensional spaces.

\section{Introduction}

The concept of extension dimension was introduced by
Dranishnikov \cite{d:94}
(see also \cite{ch1:97}, \cite{dd:96}). For a normal space
$X$ and a
space $K$ we write $e-dim X\leq K$ (the extension
dimension of $X$
does not exceed $K$) if $K$ is an absolute extensor for
$X$. This means that any continuous map $f \colon A \to K$,
defined on a closed subset $A$ of $X$, admits a continuous
extension $\bar{f} \colon X \to K$.
We can enlarge the class of normal spaces $X$ with $e-dim X\leq K$ by introducing
the following notion (see \cite[Definition 2.5]{lr:91}: for a space $K$ let
$\alpha (K)$ be the class of all normal spaces $X$ such that if $A\subset X$
is closed and
$f: A\rightarrow K$ is a continuous map, then $f$ can be extended to a map   
from $X$ into $K$ provided there is a neighborhood $U$ of $A$ and a map
$g: U\rightarrow K$ extending $f$. Obviously, $\alpha (K)$ contains all spaces
$X$ with $e-dim X\leq K$ and if $K$ is an absolute neighborhood extensor for $X$
(br., $K\in ANE(X)$), then $X\in\alpha (K)$ is equivalent to $e-dim X\leq K$
(for example, this is true if $X$ is metrizable, or more general when $X$
admits a perfect map onto a first countable space \cite{jd:95}). 

A space $X$ is said to be $\mathcal P$-like \cite{du:97}, where $\mathcal P$
is a given class, if for every open locally finite cover $\omega$ of $X$ there exists an
$\omega$-map from $X$ into a space $Y\in\mathcal P$. It is well known that a
normal space $X$ satisfies $dimX\leq n$ if and only if
$X$ is $\mathcal P$-like with $\mathcal P$ being the class of all simplicial
complexes of dimension $\leq n$. It follows from \cite{du:97} and \cite{lrs:98}
that if $X$ is compact and $K$ is a $CW$-complex, then $e-dim X\leq K$ iff $X$ is
$\mathcal P$-like, where $\mathcal P$ denotes the class of compact metric spaces
$Y$ with $e-dim Y\leq K$.

In the present note we prove that for a $CW$-complex $K$ we have $X\in\alpha (K)$
provided $X$ is $\mathcal P$-like with respect to the class of paracompact
spaces from $\alpha (K)$. In particular, if a metrizable space $X$ is
$\mathcal P$-like
with respect to metrizable spaces $Y$ with $e-dim Y\leq K$, then $e-dim X\leq K$.
It is also shown that, for countable complexes $K$, this property characterizes
the class $\alpha (K)$. We apply the above results to show that  
if $f:X\rightarrow Y$ is a refinable map between
metrizable spaces, then $e-dim X\leq K$ iff $e-dim Y\leq K$, $K$ is any $CW$-complex.
The last generalizes A. Koyama's results from \cite{ak:94}, as well as a
result of A. Koyama and R. Sher \cite{ks:94}. In the final section we prove
that refinable maps preserve $S$-weakly infinite-dimensionality (resp., finite
$C$-space property). In the class of compact spaces this result was 
proved earlier by H. Kato \cite{k:83} and A. Koyma \cite{ak:84}
(resp., D. Garity and D. Rohm \cite{gr:86}).

All spaces in this paper are assumed to be normal and all maps are continuous.

\section{The class $\alpha (K)$}

Let start with the following lemma.

\medskip
{\sc Lemma 2.1.} {\em Let $L$ be a space homotopy equivalent to a space $K$. Then 
$X\in\alpha (L)$ if and only if $X\in\alpha (K)$.} 

\medskip
{\sf Proof.}
There exist maps $\varphi \colon K\to L$ and $\psi \colon L\to K$
such that $\varphi\circ\psi\simeq id_L$ and $\psi\circ\varphi\simeq id_K$.
It suffices to show that $X\in\alpha (L)$ yields $X\in\alpha (K)$. To this
end, let $f\colon A\to K$ be a map such that $A\subset X$ is closed and $f$
can be extended to a map $g\colon U\to K$, where $U$ is a neighborhood of $A$
in $X$. Since $X$ is normal, there exists open $V\subset X$ with
$A\subset V\subset\overline{V}\subset U$. Then $k=(\varphi\circ g)|\overline{V}$ is a map
from $\overline{V}$ into $L$ which is extendable to a map from $U$ into $L$.
Using that $X\in\alpha (L)$ we can find an extension
$h\colon X\to L$ of $k$.
Because $(\psi\circ h)|V\simeq g$, by \cite[IV, 2.1]{hu:65}, $f$ admits an
extension. Hence $X\in\alpha (K)$. $\Box$

\medskip 
Recall that a map $f: X\rightarrow Y$, where $X$ and $Y$ are
metrizable spaces, is called 
uniformly 0-dimensional \cite{k:52} if there exists a compatible
metric on $X$ such that 
for every $\varepsilon >0$ and every $y\in f(X)$ there is an open
neighborhood $U$ of $y$ such
that $f^{-1}(U)$ can be represented as the union of disjoint
open sets of $diam <\varepsilon$.
It is well known that every metric space admits a uniformly
0-dimensional map into Hilbert 
cube $Q$.  

\medskip
{\sc Theorem 2.2.}
{\em Let $K$ be a countable $CW$-complex and $f \colon X\to Y$, where
$X\in\alpha (K)$ and $Y$ is metrizable. Then there exist metrizable $Z$ and
maps $h \colon X\to Z$ and $g \colon  Z\to Y$ such that $e-dim Z\leq K$ and
$f=g\circ h$.}

\medskip
{\sf Proof.}
Let $k\colon Y\to Q$ be a uniformly 0-dimensional map. Fix a Polish $ANR$-space
$P$
homotopy equivalent to $K$. By Lemma 2.1, $X\in\alpha (P)$. Since $P\in ANE(X)$,
we have $e-dim X\leq P$. Then, according to \cite[Proposition 4.9]{ch1:97},
there exist a separable metric space $M$ with $e-dim M\leq P$ and
maps $\psi \colon X\to M$, $\varphi \colon M\to Q$ such that
$\varphi\circ\psi =k\circ f$. Applying once more Lemma 2.1, we conclude that
$e-dim M\leq K$.
Let $Z$ be the fibered product of $M$ and $Y$ with respect to $\varphi$ and
$k$, and let $q\colon Z\to M$ and $g\colon Z\to Y$ denote the corresponding
projections. Since $k$ is uniformly 0-dimensional, we can find a compatible
metric $d$ on $Z$ such that $q$ is uniformly 0-dimensional with respect to
$d$ (see \cite{ap:73}). Because $e-dim M\leq K$, the last yields $e-dim Z\leq K$
\cite[Theorem 1]{l:98}. Finally, define $h\colon X\to Z$ by
$h(x)=(\psi (x),f(x))$. $\Box$

\medskip  
A map $f \colon X\to Y$ is called an $\omega$-map, where $\omega$ is an
open cover of $X$, if there is an open cover $\gamma$ of $Y$ such that
$f^{-1}(\gamma)$ refines $\omega$.

\medskip
{\sc Corollary 2.3.}
{\em Let $K$ be a countable $CW$-complex and $X\in\alpha (K)$. Then for
every locally
finite open cover $\omega$ of $X$ there is an $\omega$-map from $X$ onto a
metrizable space $Z$ with $e-dim Z\leq K$.}

\medskip
{\sf Proof.}
Take an $\omega$-map $f$ from $X$ into a metrizable space $Y$ and apply
Theorem 2.2 to obtain a metrizable space $Z$ and maps $h\colon X\to Z$,
$g\colon Z\to Y$ such that  $e-dim Z\leq K$ and $f=g\circ h$. It remains to
observe that $h$ is an $\omega$-map because $f$ is such a map. $\Box$

\medskip
It is certainly true that if Theorem 2.2 holds for any $CW$-complex, then
Corollary 2.3 also holds for arbitrary $K$.

\medskip
{\sc Theorem 2.4.}
{\em Let $X$ be a normal space and $K$ be a $CW$-complex. If for every
locally
finite open cover $\omega$ of $X$ there exists an $\omega$-map into a
paracompact space $Y$ with $Y\in\alpha (K)$, then $X\in\alpha (K)$.}

\medskip
{\sf Proof.}
We follow the ideas from the proof of \cite[Lemma 2.1]{du:97}.
There exists a normed space $Z$ and an open subset $L$ of $Z$ homotopy
equivalent to $K$.
According to Lemma 2.1, it suffices to show that $X\in\alpha (L)$.
To this end, take a map $f\colon A\to L$, where $A\subset X$ is closed,
such that $f$ is extendable to a map from a neighborhood $U$ of $A$ into
$L$. Since $X$ is normal, we can suppose that $f$ is a map from $U$ into $L$
and there is an open set $V\subset X$ containing $A$ such that
$\overline{V}\subset U$. To prove that $f|A$ can be extended to a map from $X$
into $L$, it suffices to find a map $\bar{f}\colon X\to L$ such that
$\bar{f}|V\simeq f|V$ (see \cite[IV, 2.1]{hu:65}). Let $\lambda$ be a locally
finite open cover of $L$ such that $St(z,\lambda)$ is contained in an open ball
for every $z\in L$. Then
$\omega =\{f^{-1}(W): W\in\lambda\}\cup\{X-\overline{V}\}$ is a locally finite
open cover of $X$, so there exists an $\omega$ map $p\colon X\to Y$ into a
paracompact space $Y$ with $Y\in\alpha (K)$. Let $\gamma$ be a locally finite
open cover of $Y$ such that $p^{-1}(\gamma)$ refines $\omega$, and
$\{s_G:G\in\gamma\}$ be a partition of unity subordinated to $\gamma$.
Consider the closure $F$ of $p(V)$ in $Y$ and let
$\gamma _F=\{G\in\gamma :G\cap F\ne\emptyset \}$. For every $G\in\gamma _F$
fix a point $x_G\in p^{-1}(G)\cap V$.
Denote $a_G=f(x_G)$ if $G\in\gamma _F$ and $a_G=0$ otherwise, where 0 is the
zero-point of $Z$.

\medskip
{\sc Claim.} {\em For every $y\in F$ there is a ball $B\subset L$ 
with $f(p^{-1}St(y,\gamma))\subset B$.}

\medskip
If $St(y,\gamma)=\{G_i: i=1,2,..,n\}$ and $G_y=\cap\{G_i:i=1,..,n\}$, then
$G_y$ meets $p(V)$, so we can find $x\in p^{-1}(G_y)\cap V$. Since each
$p^{-1}(G_i)$ is contained in a element of $\omega$ and meets $V$, we have
$p^{-1}(St(y,\gamma)\subset St(x,f^{-1}(\lambda))$. Hence 
$f(p^{-1}(St(y,\gamma))\subset St(f(x),\lambda)$. Finally, choose an open
ball $B\subset L$ with $St(f(x),\lambda)\subset B$.

\medskip
Define a map $g\colon F\to L$ by $g(y)=\sum\{s_G(y)\cdot a_G:G\in\gamma _F\}$.
According to the Claim, this definition is correct. Moreover, again by the
Claim, for every $x\in V$ there is an open ball in $L$ containing both $f(x)$
and $g(p(x))$. Therefore $f|V\simeq (g\circ p)|V$. Observe that the formula
$\bar{g}=\sum\{s_{G}(y)\cdot a_G: G\in\gamma\}$ defines a map from $Y$ into $Z$
which extends $g$. Obviously $O=\bar{g}^{-1}(L)$ is an open subset of $Y$
and $g$ has an extension to a map from $O$ into $L$.
Since $Y\in\alpha (K)$, by Lemma 2.1, $Y\in\alpha (L)$.
Hence $g$ can be extended to a map $q\colon Y\to L$. Let $\bar{f}=q\circ p$.
Then $\bar{f}|V=(q\circ p)|V=(g\circ p)|V\simeq f|V$. $\Box$

\medskip
Combining Corollary 2.3 and Theorem 2.4 we obtain the following 
characterization of the class $\alpha (K)$.

\medskip
{\sc Corollary 2.5.}
{\em Let $K$ be a countable $CW$-complex. Then a normal space $X$ belongs to
$\alpha (K)$ if and only if for every locally finite open cover $\omega$ of
$X$ there exists an $\omega$-map into a metrizable space $Y$ with
$e-dim Y\leq K$.}

\section{Extension dimension}

A surjective map $r\colon X\to Y$ is called refinable \cite{ak:94} if for any
locally finite 
open cover $\omega$ of $X$ and  any open cover $\gamma$ of $Y$ there exists a
surjective $\omega$-map $f\colon X\to Y$ such that $r$ and $f$ are
$\gamma$-close (i.e., for every $x\in X$ there is an element of $\gamma$
containing both points $r(x)$ and $f(x)$). The map $f$ is called
$(\omega ,\gamma )$-refinement of $r$. When there exists a closed
$(\omega ,\gamma )$-refinement for $r$ we say that $r$ is  c-refinable.
For compact metric spaces this definition coincides with the original
one given by J. Ford and J. Rogers \cite{fr:78}

Koyama proved that $dimX=dimY$ provided $X$ and $Y$ are metric spaces and $r$
is refinable \cite[Theorem 2]{ak:94}.
Under the same hypotheses A. Koyama and R. Sher showed \cite{ks:94}
that $dim_GX=dim_GY$ for any finitely generated Abelian group $G$. Here
$dim_GX$ stands for the cohomological dimension of $X$ with respect to $G$.
In case $r$ is c-refinable we have $K\in AE(X)$ if and only if $K\in AE(Y)$
for any simplicial complex $K$ (see \cite[Theorem 1]{ak:94}).
In the present section we generalize all these results by proving that if
$r$ is a refinable map between metric spaces $X$ and $Y$, then $e-dim X\leq K$
if and only if $e-dim Y\leq K$ for any $CW$-complex $K$.

\medskip
{\sc Lemma 3.1.}
{\em Let $r$ be a refinable map from a normal space $X$ onto a paracompact
space
$Y$ and $L$ a locally convex subset of a normed space $Z$. Suppose
$A\subset W\subset\overline{W}\subset W_1\subset\overline{W_1}\subset F\subset W_2$,
where $W$, $W_1$ and $W_2$ $($resp., $A$ and $F$$)$ are open $($resp., closed$)$
subsets of $Y$. Further,
let $H_1=r^{-1}(\overline{W_1})$ and $H_2=r^{-1}(\overline{W_2})$. If
$g\colon\overline{W_1}\to L$
and $h\colon H_2\to L$ are continuous with $h|H_1=(g\circ r)|H_1$, then $g|A$
can be extended to a continuous map from $F$ into $L$.}

\medskip
{\sf Proof.}
Obviously, $L$ is an $ANR$ as a set having a base of convex subsets. By
\cite[IV, 2.1]{hu:65}, $g|A$ is extendable over $F$ provided there exists a map
$p\colon F\to L$ such that $p|D\simeq g|D$ for some neighborhood $D$ of $A$.
We shall construct a map $p$ satisfying this condition with $D=W$.  
 
Let $\lambda$ be a locally finite
open cover of $L$ such that $St(G,\lambda)$ is contained in a
convex subset of $L$ for
every $G\in\lambda$. Take an open set $U\subset Y$ such that
$F\subset U\subset\overline{U}\subset W_2$ and denote
$\omega =\{h^{-1}(G)\cap r^{-1}(W_2): G\in\lambda\}\cup
\{X-r^{-1}(\overline{U})\}$.
Take a locally finite open cover $\gamma$ of $Y$ refining the cover
$\{g^{-1}(G)\cap W_1:G\in\lambda\}\cup\{Y-\overline{W}\}$ such that
$St(\overline{W},\gamma )\subset W_1$ and $St(F,\gamma )\subset U$.
There exists a $(\omega ,\gamma )$-refinement
$f\colon X\to Y$ for $r$, and let $\beta$ be an
open star-refinement of $\gamma$ such that $f^{-1}(\beta )$ refines $\omega$.
Since $f$ is surjective, each $f^{-1}(V)\ne\emptyset$, $V\in\beta$.
Let $\beta _1=\{V\in\beta :V\cap F\ne\emptyset\}$ and for every
$V\in\beta _1$ pick a point $x_V\in f^{-1}(V\cap F)$ such that
$x_V\in f^{-1}(W)$ when $V$ meets $W$.

\medskip
{\sc Claim.} {\em $f^{-1}(W)\subset H_1$ and $f^{-1}(F)\subset r^{-1}(U)$}.

If $z\in f^{-1}(W)$, then there exists $O\in\gamma$ such that $f(z),r(z)\in O$
(recall that $f$ is $\gamma$-close to $r$). So, $O\cap W\neq\emptyset$, i.e.
$O\subset St(\overline{W},\gamma)$. Since
$St(\overline{W},\gamma)\subset W_1$,
$O\subset W_1$. Hence, $z\in r^{-1}(O)\subset H_1$. Similarly, using the
inclusion $St(F,\gamma )\subset U$, we can show that
$f^{-1}(F)\subset r^{-1}(U)$.

\medskip
Since, $f^{-1}(F)\subset r^{-1}(U)$ (see the claim above), $f^{-1}(V)$ meets
$r^{-1}(U)$ for every $V\in\beta _1$. On the other hand, $f^{-1}(V)$ is
contained in an element of $\omega$, hence each $f^{-1}(V)$, $V\in\beta _1$,
is contained in $r^{-1}(W_2)$. Therefore, the points 
$a_V=h(x_V)$, $V\in\beta _1$, are determined. Consider also the points
$y_V=f(x_V)$ for $V\in\beta _1$. Take a
partition of unity $\{s_V:V\in\beta\}$ subordinated to $\beta$ and define the
map $p\colon F\to Z$ by $p(y)=\sum\{s_V(y)\cdot a_V:V\in\beta _1\}$.

\medskip
Let show that $p$ is a map from $F$ into $L$. If $y\in F$ and $x\in X$ with
$y=f(x)$, then 
$f^{-1}(St(y,\beta ))\subset St(x,\omega)\subset r^{-1}(W_2)$. So,
$h(f^{-1}(St(y,\beta ))\subset St(h(x),\lambda)$.  Consequently,
$St(h(x),\lambda)$ contains all points $a_V$ with $y\in V$. Since, there
exists 
a convex set $B\subset L$ containing $St(h(x),\lambda)$, we obtain that
$p(y)\in B$.

It remains to prove that $p|W\simeq g|W$. Since $\beta$ is a
star-refinement of $\gamma$, for any $y\in W$ there exists
$G_y\in\lambda$ with $St(y,\beta )\subset g^{-1}(G_y)$. Fix $y\in W$ and
$x\in X$ with $y=f(x)$. Let $\{V(1), V(2),..,V(n)\}$ be the set of all
$V\in\beta _1$ containing $y$. According to the choice of the points $x_V$,
each $x_{V(i)}\in f^{-1}(W)$, so
$y$ and $y_{V(i)}$, $i=1,..,n$,
belong to $W$. Now, since $f,r$ are $\gamma$-close, we can find
$G_i\in\lambda$, $i=1,..,n$ such that
$y_{V(i)},r(x_{V(i)})\in g^{-1}(G_i)$.
Therefore $y\in g^{-1}(G_y)$ and
$y_{V(i)}\in g^{-1}(G_y\cap G_i)$. Hence
$\{g(y), g(r(x_{V(i)})):i=1,..,n\}$ is a subset of $St(G_y,\lambda )$.
The last set is contained in a convex set $B_1\subset L$.
Because 
$x_{V(i)}\in f^{-1}(W)\subset H_1$, $i=1,..,n$, we have
$g(r(x_{V(i)}))=h(x_{V(i)})=a_{V(i)}$, $i=1,..,n$. 
We finally obtain that $B_1$ contains $g(y)$ and all $a_V$ with $y\in V$. 
Therefore $B_1$ contains both $p(y)$ and $g(y)$. So, $p|W\simeq g|W$. $\Box$

\medskip
{\sc Proposition 3.2.}
{\em Let $r$ be a refinable map from a normal space $X$ onto a paracompact space
$Y$. If $K$ is any $CW$-complex and $X\in\alpha (K)$, then $Y\in\alpha (K)$.}

\medskip
{\sf Proof.}
As in the proof of Theorem 2.4, let $L$ be an open subset of a normed space $Z$
homotopy equivalent to $K$. It suffices to show that $Y\in\alpha (L)$.
Towards this end, let $g\colon A\to L$ be a map with $A\subset Y$ closed and
such that $g$ is extendable to a map $\bar{g} \colon U\to L$, where $U$ is a
neighborhood of $A$. Choose open in $Y$ sets $W$ and $W_1$ such that
$A\subset W\subset\overline{W}\subset W_1\subset\overline{W_1}\subset U$ and
let $H_1=r^{-1}(\overline{W_1})$. 
Since, by Lemma 2.1, $X\in\alpha (L)$, and the map $(\bar{g}\circ r)|H_1$ is
extendable to a map from a neighborhood of $H_1$ into $L$ (the map
$\bar{g}\circ r\colon r^{-1}(U)\to L$ can serve as such an extension), there
exists a map $h\colon X\to L$ extending $(\bar{g}\circ r)|H_1$. Now, we 
apply Lemma 3.1 (with $F$ replaced by $Y$ and $H_2$ by $X$) to conclude
that $g$ is extendable to a map from $Y$ into $L$. Hence, $Y\in\alpha (L)$.
$\Box$ 

\medskip
{\sc Proposition 3.3.}
{\em Let $K$ be a $CW$-complex and $r$ be a refinable map from a normal space $X$
onto a paracompact space $Y$ with $Y\in\alpha (K)$. Then $X\in\alpha (K)$.}

\medskip
{\sf Proof.}
Since for every locally finite open cover $\omega$ of $X$ there exists an
$\omega$-map from $X$ onto $Y$, the proof follows from Theorem 2.4. $\Box$

\medskip
Proposition 3.2 and Proposition 3.3 imply the following general result. 

\medskip
{\sc Corollary 3.4.}
{\em For every $CW$-complex $K$ and a refinable map from a normal space
$X$ onto a
paracompact space $Y$ we have $X\in\alpha (K)$ if and only if
$Y\in\alpha (K)$.}

\medskip
Since every $CW$-complex is an $ANE$ for the class of all metrizable spaces,
we have

\medskip
{\sc Corollary 3.5.}
{\em If $r$ is a refinable map between the metric spaces $X$ and $Y$, and
$K$ is a
$CW$-complex, then $e-dim X\leq K$ if and only if $e-dim Y\leq K$.}

\section{Infinite-dimensional spaces}

The preservation of infinite-dimensional properties under refinable maps
is widely treated by different authors. H. Kato \cite{k:83} has shown
that refinable maps between compact metric spaces preserve weakly
infinite-dimensionality and A. Koyama \cite{ak:84} extended this result 
by proving that $S$-weakly infinite-dimensionality is preserved by
$c$-refinable maps between normal spaces. The analogous question
concerning property $C$ was settled by 
D. Garity and D. Rohm \cite{gr:86} for compact metric spaces. 
F. Ancel \cite{a:85}
introduced approximately invertible maps and established that any 
such a map with compact fibres and metric domain and range preserves
property $C$. 
Because every refinable map between compact metric spaces is 
approximately invertible, Ancel's result implies that one of 
Garity and Rohm.  
We shall prove in this section that refinable maps with normal 
domain and paracompact range preserve $S$-weakly 
infinite-dimensionality and finite $C$-property. Let us
note that in the class of compact spaces
weakly infinite-dimensionality and $S$-weakly infinite-
dimensionality, as well as $C$-space property and finite $C$-space
property, are equivalent. 

Recall that a space $X$ is called
$A$-weakly infinite-dimensional \cite{ap:73} (resp. $S$-weakly 
infinite-dimensional) if
for any sequence $\{A_i,B_i\}$ of disjoint pairs of closed sets in $X$
there exist closed separators $C_i$ between $A_i$ and $B_i$ such that
$\cap _{i=1}^{\infty}C_i=\emptyset$ (resp., 
$\cap _{i=1}^nC_i=\emptyset$ for some $n$). Usually, $A$-weakly 
infinite-dimensional spaces are called weakly infinite-dimensional. 

Another type infinite-dimensional property is the following one:
$X$ is said to be a $C$-space \cite{ag:78} if for any sequence of open
covers $\{\omega _i\}$ of $X$ there exists a sequence of disjoint open
families $\{\gamma _i\}$ such that $\gamma _i$ refines $\omega _i$ and
$\cup _{i=1}^{\infty}\gamma _i$ covers $X$. The sequence $\{\gamma _i\}$
is called a $C$-refinement of $\{\omega _i\}$.  
If, in the above definition,
the sequence $\{\gamma _i\}$ is finite (and satisfying all other
properties), then $X$ is called a finite $C$-space \cite{b1:0}. Every
finite $C$-space is $S$-weakly infinite-dimensional and every $C$-space
is weakly infinite-dimensional.
Moreover, $X$ is $S$-weakly infinite dimensional
(resp., a finite $C$-space) iff $\beta X$ is weakly infinite-dimensional
(resp., a $C$-space), see \cite{ap:73} and \cite{vc:0}. Recall
that every countable dimensional (a countable union of finite-dimensional
subspaces) metric space has property $C$, but there exists a metrizable
$C$-compactum which is not countable dimensional \cite{rp:81}. 

Let us note that, even for the class of compact metric spaces, 
there is no $CW$-complex $K$ such that $X$ is weakly infinite-
dimensional (resp., $X$ has property $C$) if and only if $e-dim X\leq K$.
Otherwise, since the Hilbert cube $Q$ is the inverse limit of
a sequence of finite-dimensional spaces, by \cite{rs:98}, we would
have that $Q$ is weakly infinite-dimensional (resp., $C$-space). 
Therefore, we can not apply the results in Section 2 to conclude that
weakly infinite-dimensionality and the property $C$ are preserved
by refinable maps. 

\medskip
{\sc Theorem 4.1.}
{\em $S$-weakly infinite-dimensionality is preserved by refinable maps with
normal domains and paracompact ranges.}

\medskip
{\sf Proof.}
Let $r$ be a refinable map from $X$ onto $Y$, where $X$ is normal
$S$-weakly infinite-dimensional and
$Y$ paracompact. We need the following characterization of $S$-weakly
infinite-dimensionality \cite{ap:73}: a space $Z$ is $S$-weakly infinite-
dimensional iff for every map $f\colon Z\to Q$ there is $n$ such that
$\pi _n\circ f$ is an inessential map; here each 
$\pi _n\colon Q\to I^n$ is the projection from $Q$ onto its $n$-dimensional
face $I^n$ generated by first $n$ coordinates.
So, fix $f\colon Y\to Q$. Then, by the above characterization, there
exists $n$ with $h_1=\pi _n\circ f\circ r$ inessential. Let $f_n=\pi _n\circ f$
and $A=f_n^{-1}(S^{n-1})$, $S^{n-1}$ being the boundary of $I^n$.
We are going to show that $f_n$ is inessential, i.e.
$f_n|A$ can be extended to a map from $Y$ into $S^{n-1}$.

To this end, fix an interior point $a$ from $I^n$ and let 
$U=f_n^{-1}(L)$ with $L=I^n-\{a\}$.
Take open sets $W$ and $W_1$ in $Y$ with
$A\subset W\subset\overline{W}\subset W_1\subset\overline{W_1}\subset U$
and denote $H_1=r^{-1}(\overline{W_1})$. Since $h_1$ is inessential and
$h_1(H_1)\subset L$, there exists a map $h\colon X\to L$ extending $h_1|H_1$.
By Lemma 3.1, $f_n|A$ admits an extension $g\colon Y\to L$, and that 
suffices to find an extension of $f_n|A$ from $Y$ into $S^{n-1}$. $\Box$

\medskip
Next proposition extends the result of D. Garity and D. Rohm \cite{gr:86}
that property $C$ is preserved by refinable maps between compact metric
spaces.

\medskip
{\sc Proposition 4.2.}
{\em Let $r\colon X\to Y$ be refinable with $X$ normal and $Y$ paracompact.
If $X$ is a finite $C$-space, then $r(K)$ has property $C$ for every
compact $K\subset X$.}

\medskip
{\sf Proof.}
Let $exp(Q)$ be the space of all closed subsets of $Q$ with the Vietoris
topology and ${\mathcal Z}(Q)\subset exp(Q)$ consist of all $Z$-sets in $Q$.
We need the following result of Uspenskij
\cite[Theorem 1.4]{vu:98}: A compact space $Z$ has property $C$ if and
only if for any
map $\phi \colon Z\to{\mathcal Z}(Q)$ there exists a map
$g\colon Z\to Q$ such that $g(z)\not\in\phi (z)$ for every $z\in Z$.

So, fix a map $\phi\colon r(K)\to {\mathcal Z}(Q)$ and take an extension
$\Phi\colon Y\to exp(Q)$ of $\phi$ (such an extension exists because
$exp(Q)$ is an $AR$). Next,
let $X_1=r_1^{-1}(Y)$ and $K_1=r_1^{-1}(r(K))$, where 
$r_1=\beta r\colon \beta X\to \beta Y$. Since $\beta X$ is a $C$-space, 
so is $K_1$ (as a closed subset of $\beta X$). 
Consider the map $\Psi\colon X_1\to exp(Q)$, $\Psi =\Phi\circ r_1$.
Then, by mentioned above result of Uspenskij, there exists a map
$h\colon K_1\to Q$ with $h(x)\not\in\Psi (x)$ for every $x\in K_1$.
Take $\epsilon >0$ such that $d(h(x),\Psi (x))>\epsilon$ for all
$x\in K_1$, where $d$ is the metric in $Q$, and extend $h$ to a map
$g\colon U\to Q$, $U$ is a neighborhood of $K_1$ in $X_1$, satisfying
$d(g(x),\Psi (x))>\epsilon$ for every $x\in U$. Next step is to find
a neighborhood  $W$ of $r(K)$ in $Y$ with
$r_1^{-1}(W)\subset r_1^{-1}(\overline{W})\subset U$ and
choose an open cover $\gamma$ of $Y$ such that
$St(r(K),\gamma)\subset W$ and $d_H(\Phi (y),\Phi (z))<2^{-1}\epsilon$
for any two $\gamma$-close points $y,z\in Y$, where $d_H$ is the
Hausdorff metric on $exp(Q)$. Let $\lambda$ be a finite open cover
of $Q$
such that each $St(q,\lambda)$, $q\in Q$, is contained in a convex
set of diameter $\leq 2^{-1}\epsilon$. Then
$\omega =\{g^{-1}(G)\cap X:G\in\lambda\}\cup\{X-r_1^{-1}(\overline{W})\}$
covers $X$. So,
there exists a $(\omega,\gamma)$-refinement $f\colon X\to Y$ of $r$ and
a locally finite open cover $\alpha$ of $Y$ with $f^{-1}(\alpha)$ refining
$\omega$. Since $St(r(K),\gamma)\subset W$, we have
$f^{-1}(r(K))\subset r^{-1}(W)$ (see the claim from Lemma 3.1). Hence,
if $\alpha _1=\{V\in\alpha :V\cap r(K)\ne\emptyset\}$, then we can choose
a point $x_V\in f^{-1}(V)\cap r^{-1}(W)$ for every $V\in\alpha _1$. 
Finally, let $\{s_V:V\in\alpha\}$ be a partition of unity
subordinated to $\alpha$ and define the map $p\colon r(K)\to Q$
by $p(y)=\sum\{s_V(y)\cdot a_V:V\in\alpha _1\}$, where
$a_V=g(x_V)$ for every $V\in\alpha _1$.

It remains only to show that $p(y)\not\in\phi (y)$ for every $y\in r(K)$.
Fix $y\in r(K)$ and $x\in X$ with $y=f(x)$. Then
$f^{-1}(St(y,\alpha))\cap r^{-1}(W)$ contains $x$ and all $x_V$ with $y\in V$.
Moreover, $f^{-1}(St(y,\alpha))\cap r^{-1}(W)$ is a subset of
$St(x,\omega)\cap r^{-1}(W)$. So,
$g(f^{-1}(St(y,\alpha))\cap r^{-1}(W))\subset St(g(x),\lambda)$.
Take a convex set $B\subset Q$ of diameter $\leq 2^{-1}\epsilon$ which
contains $St(g(x),\lambda)$. Then, $p(y),g(x)\in B$ and, because
$d(g(x),\Psi (x))>\epsilon$, $d(p(y),\Psi (x))>2^{-1}\epsilon$.
On the other hand, $y=f(x)$ and $r(x)$ are $\gamma$-close, so
$d_H(\phi (y),\Psi (x))<2^{-1}\epsilon$. Therefore, $p(y)\not\in\phi (y)$. $\Box$

\medskip
{\sc Theorem 4.3.}
{\em Finite $C$-space property is preserved by refinable maps with
normal domains and paracompact ranges.}

\medskip
{\sf Proof.}
Suppose $X$ is a normal finite $C$-space, $Y$ is paracompact and
$r\colon X\to Y$ is refinable. By \cite{b1:0}, there exists a compact
$C$-space $K\subset X$ such that every closed set $H\subset X$ which is disjoint
from $K$  has a finite dimension $dim$. Because this property characterizes
finite $C$-spaces \cite{b1:0}, it suffices to show that every closed set in
$Y$ disjoint from $r(K)$ is finite-dimensional and $r(K)$ has property $C$.
By Proposition 4.2, $r(K)$ is a $C$-space. So, it remains only to show that
all closed sets in $Y$ outside $r(K)$ are finite-dimensional.

Suppose $B\subset Y$ is
closed and $B\cap r(K)=\emptyset$. Take an open $W_2\subset Y$ containing $B$
with $\overline{W_2}$ disjoint from $r(K)$. Then $H_2=r^{-1}(\overline{W_2})$
is closed in $X$ and disjoint from $K$. Hence, $dim H_2=n$ is finite. We
shall prove that $dim B\leq n$.

Let $A\subset B$ be closed and $q\colon A\to S^n$. Our intention is
to extend $q$  to a map from $B$ into $S^n$.
Take an open set
$U\subset Y$ such that $B\subset U\subset\overline{U}\subset W_2$
and denote $F=\overline{U}$. The proof will be completed if we can find
an extension $\bar{q}\colon F\to L$ of $q$, where $L$ stands for the cube
$I^{n+1}$ with deleted an interior point. Towards this end, extend $q$ to a map
$g\colon\overline{W_1}\to S^n$, where $W_1\subset Y$ is open with
$A\subset W_1\subset\overline{W_1}\subset U$, and then choose any open set
$W\subset Y$
satisfying $A\subset W\subset\overline{W}\subset W_1$. Then
$g\circ r$ is a map from $H_1=r^{-1}(\overline{W_1})$ to $S^n$ and,
since $dim H_2\leq n$, there exists an extension $h\colon H_2\to L$
of $g\circ r$. Finally, apply Lemma 3.1, to find an extension
$f\colon F\to L$ of $g$ and observe that $f|A=q$. $\Box$

\bigskip\noindent
Department of Mathematics and Statistics, 
University of Saskatchewan,

\noindent
McLean Hall,
106 Wiggins Road, Saskatoon, SK, Canada S7N 5E6

E-mail: chigogid@math.usask.ca

\medskip\noindent
Department of Mathematics,
University of Swaziland, Pr. Bag 4, 

\noindent
Kwaluseni, Swaziland 

E-mail: valov@realnet.co.sz

\end{document}